\documentclass[10pt,reqno,a4paper]{amsart}

\newcommand{\K}{k}

\newcommand{\Rb}{\mathbb{R}}
\newcommand{\SXb}{\mathbb{S}}

\newcommand{\cyb}{\operatorname{CYB}}

\newcommand{\mink}{\mathbb{M}^{(1,d-1)}}
\newcommand{\til}[1]{\widetilde{#1}}
\newcommand{\dd}{\operatorname{d}\!}
\newcommand{\dep}[2]{\frac{\partial #1}{\partial #2}}
\newcommand{\er}{\mathcal{R}}
\newcommand{\zed}{\mathcal{Z}}
\newcommand{\boite}[2]{\;\boxed{\overset{#1}{\underset{#2}{}}}\;}

\newcommand{\ush}{\operatorname{USh}}
\newcommand{\refeqn}[1]{(\ref{#1})}
\newcommand{\cp}{\operatorname{cp}}
\newcommand{\alt}{\operatorname{Alt}}
\newcommand{\ggo}{\mathfrak{g}}
\newcommand{\id}{\operatorname{id}}
\newcommand{\eps}{\varepsilon}
\newcommand{\h}{h}
\newcommand{\sfo}[1]{#1[[\h]]}
\newcommand{\bproof}{\emph{Proof. }}
\newcommand{\eproof}{\hfill\qed\bigskip}
\newcommand{\Hom}{\operatorname{Hom}}
\newcommand{\Lm}{\operatorname{L}}
\newcommand{\Rm}{\operatorname{R}}
\newcommand{\et}{\quad\text{and}\quad}
\newcommand{\bw}{\wedge}

\theoremstyle{plain}
\newtheorem{thm}{Theorem}[section]
\newtheorem{prop}[thm]{Proposition}

\newtheorem{lem}[thm]{Lemma}

\theoremstyle{definition}
\newtheorem{defi}[thm]{Definition}

\theoremstyle{remark}
\newtheorem{rmk}[thm]{Remark}

\begin{document}

\title[A QLBA formulation of the PR Poisson algebra]{A quasi-Lie bialgebra formulation\\of the Pohlmeyer-Rehren Poisson algebra}

\author[M. Bordemann]{Martin Bordemann}
\address{Martin Bordemann: LMIA, Universit\'{e} de Haute-Alsace, 4 rue des Fr\`{e}res Lumi\`{e}re, F-68093 Mulhouse}
\email{martin.bordemann@uha.fr}

\author[B. Enriquez]{Benjamin Enriquez}
\address{Benjamin Enriquez: IRMA (CNRS), Universit\'e Louis Pasteur, 7 rue Ren\'e Descartes, F-67084 Strasbourg}
\email{enriquez@math.u-strasbg.fr}

\author[L. Hofer]{Laurent Hofer}
\address{Laurent Hofer: LMIA, Universit\'{e} de Haute-Alsace, 4 rue des Fr\`{e}res Lumi\`{e}re, F-68093 Mulhouse}
\email{laurent.hofer@uha.fr}

\begin{abstract}
We present a quasi-Lie bialgebra (QLBA) quantization problem
which comes from an algebraic reformulation of the Nambu-Goto
string theory and invariant charges by Pohlmeyer and Rehren.
This QLBA structure depends on a symmetric
bivector (coming from a Minkowski metric) and is built on the
free Lie algebra on a finite dimensional vector space. We solve
this problem when the bivector has rank 1 or 2.
\end{abstract}

\maketitle
\thispagestyle{empty}

\section*{Introduction}

The quantization of a quasi-Lie bialgebra $(\ggo,\delta,\varphi)$
(see \cite{Dri90}) is an open problem, however if the $\varphi=0$
(\emph{i.e.} $(\ggo,\delta)$ is a Lie bialgebra) we can use the deep
result of Etingof-Kazhdan \cite{EK96,EK98a,EK98b}. The aim of this
paper is to present a quasi-Lie bialgebra quantization problem which
comes from the theory of reparametrization-invariant charges due to
Pohlmeyer-Rehren \cite{PR86}. Details of the constructions presented
here can also be found in the third author's PhD-thesis \cite{Hof07}.
\begin{itemize}
\item In Section \ref{sec:2}, we recall the basic material for
    quasi-Lie bialgebras (QLBAs), quasi-Hopf algebras and their deformation quantization.
    We emphasize the r\^{o}le of traces on the universal envelopping algebra
    of a quasi Lie bialgebra and on a quasi-Hopf algebra, respectively: this space
    becomes a Poisson algebra and an associative algebra, respectively.
\item In Section \ref{sec:3}, we introduce a QLBA
   structure on the free Lie algebra over a vector space $V$ related to a bivector
   in $V\otimes V$. In case the bivector has rank $1$ or $2$
   this QLBA structure turns out to be a bialgebra
   structure and can explicitly be quantized by a simple exponential formula.
   For higher rank we compute a solution of the quantization problem
   modulo $\h^3$.
\item Finally, in Section \ref{sec:1}, we briefly sketch the theory of invariant
   charges $\zed_{\mu_1\cdots\mu_k}^\pm$ of closed Nambu-Goto strings due to Pohlmeyer,
   and the Poisson structure of these quantities due to Pohlmeyer and
   Rehren. Moreover we establish the coincidence of their Poisson structure with
   the one we have established in the framework of QLBAs in Section
   \ref{sec:3}.
   The link with a quasi-co-Poisson structure is established in
   Proposition \ref{PPoRehQuasi}.
\end{itemize}

\section{Algebraic framework of quasi-Lie bialgebras and their quantization}\label{sec:2}

The material presented here is standard, see \cite{Dri90} for further details.
In this section $\K$ is a  field of characteristic 0 and $\sfo\K$ the ring of
formal series on $\K$. In the following, $\cp$ denotes the sum over all
cyclic permutations in three arguments, and $\alt$ denotes the alternating
sum over all permutations.

\subsection{QLBA and their associated quasi-(co)Poisson algebras}\label{rmk:quasi_co_poisson}

According to Drinfel'd, a quasi-Lie bialgebra (QLBA in the sequel)
is a Lie algebra $(\ggo,[,])$ equipped with a $1$-cocyle
$\delta:\ggo\rightarrow\bw^2\ggo$ such that the $1$-cocycle
$\cp(\delta\otimes\id)\delta:\ggo\rightarrow\bw^3\ggo$ is the coboundary of
an element $\varphi\in\bw^3\ggo$:
\begin{equation}
 \cp(\delta\otimes\id)\delta(x)
 = [x\otimes 1\otimes 1+1\otimes x\otimes 1+1\otimes 1\otimes x,\varphi]\label{eqn:cojac}
\end{equation}
with $x\in\ggo$.
Moreover the additional condition
$\alt(\delta\otimes\id\otimes\id)(\varphi)=0$ is required.
We denote such a structure as a triple $(\ggo,\delta,\varphi)$.

Let $U(\ggo)$ be the universal envelopping algebra of the Lie algebra
$\ggo$. It is known to be a Hopf algebra with the cocommutative
shuffle comultiplication
$\Delta_0:U(\ggo)\rightarrow U(\ggo)\otimes U(\ggo)$ induced by
$\Delta_0(x)=x\otimes 1+1\otimes x$ for all $x\in\ggo$.
 It is well-known that the $1$-cocycle $\delta$ can be extended from $\ggo$ to $U(\ggo)$
as a derivation $D:U(\ggo)\rightarrow U(\ggo)\otimes U(\ggo)$ along
the algebra homomorphism $\Delta_0$ (see subsection \ref{sec:der} for definitions).
The values of $D$ are in $\wedge^2 U(\ggo)$.
It satisfies the \textit{co-Leibniz rule}:
\begin{equation}\label{eqn:coleibniz}
 (\Delta_0\otimes\id)D = (\id\otimes D)\Delta_0+\sigma_{23}(D\otimes\id)\Delta_0
\end{equation}
where $\sigma_{23}(x_1 \otimes x_2 \otimes x_3) = x_1 \otimes x_3\otimes x_2$,
and the \textit{quasi-co-Jacobi identity}:
\begin{equation}\label{eqn:QuasiJacobi}
 \cp(D\otimes\id)D(x) = [(\id\otimes\Delta_0)\Delta_0(x),\varphi]
\end{equation}
for all $x\in U(\ggo)$.
Such a $D:U(\ggo)\rightarrow \bw^2 U(\ggo)$ is called a \textit{quasi-co-Poisson
bracket}, and $(U(\ggo),\Delta_0,D)$ is called a \textit{quasi-co-Poisson
algebra}.

Let $U(\ggo)^*$ the algebraic dual of $U(\ggo)$. For $F,G\in U(\ggo)^*$
there are two operations, i.e.
\begin{equation}\label{eqn:MultShuffleUGBracketD}
   F\bullet_0 G:= (F\otimes G) \Delta_0~~~\mathrm{and}~~~
   \{F,G\}_D:=(F\otimes G) D
\end{equation}
where $\bullet_0$ defines an associative
commutative mutiplication, whereas the bracket $\{\,,\}_D$ is bilinear
and antisymmetric. The co-Leibniz rule above implies
that $\{\,,\}_D$ satisfies the usual \textit{Leibniz rule}. The Jacobi identity does
\textit{not} hold in general.\\
However, let
$\operatorname{Tr}\big(U(\ggo)\big)\subset U(\ggo)^*$ be the subspace of
all \textit{traces}, i.e.
\begin{equation}\label{eqn:DefTraces}
   \operatorname{Tr}\big(U(\ggo)\big)
   =\{F\in U(\ggo)^*~|~ F(ab)=F(ba)~~ \forall~a,b\in U(\ggo)\}.
\end{equation}
It is an easy consequence of the
fact that $\Delta_0$ is a homomorphism of associative algebras that
$\operatorname{Tr}\big(U(\ggo)\big)$ is stable under $\bullet_0$ and
$\{\,,\}_D$. Moreover eqn (\ref{eqn:QuasiJacobi}) implies
that the Jacobi identity for the bracket $\{\,,\}_D$ restricted to
$\operatorname{Tr}\big(U(\ggo)\big)$ holds, whence
\begin{prop}\label{PropPoissonOnTraces}
$\big(\operatorname{Tr}\big(U(\ggo)\big),\bullet_0,\{\,,\}_D\big)$ is a
Poisson algebra called the \textit{associated Poisson algebra to the QLBA}
$(\ggo,[\,,],\delta,\varphi)$.
\end{prop}

\subsection{Quasi-Hopf algebras}
A quasi-Hopf algebra (QH algebra in the sequel) is an algebra $(A,\mu,\eta)$
equipped with a coproduct $\Delta:A\rightarrow A\otimes A$ and a counit
$\eps:A\rightarrow\K$ such that the two $\K$-morphisms
$(\id\otimes\Delta)\Delta$ and $(\Delta\otimes\id)\Delta$ are
conjugated by an invertible element $\Phi\in A^{\otimes 3}$:
\begin{equation}
 (\id\otimes\Delta)\Delta(x)\Phi = \Phi(\Delta\otimes\id)\Delta(x)\label{eqn:coass}
\end{equation}
with $x\in A$. There are also some additional conditions
$(\eps\otimes\id)\Delta=(\id\otimes\eps)\Delta=\id$,
$(\id\otimes\eps\otimes\id)(\Phi)=1\otimes 1$ and the so called pentagon condition:
\begin{equation*}
 \Phi^{1,2,34}\Phi^{12,3,4}=\Phi^{2,3,4}\Phi^{1,23,4}\Phi^{1,2,3}
\end{equation*}
where $\Phi^{1,2,34}$, $\Phi^{2,3,4}$,... have the obvious meanings
$(\id\otimes\id\otimes\Delta)(\Phi)$, $1\otimes\Phi$,... Moreover
the existence of a quasi-antipode $(S,\alpha,\beta)$ is required
(see (1.17), (1.18) and (1.19) in \cite{Dri90}). We denote such a
structure as a quadruple $(A,\mu,\Delta,\Phi)$.

Again, let $A^*$ be the algebraic dual of $A$, and let
$\operatorname{Tr}(A)\subset A^*$ be the set of traces of $A$, i.e.
\[
  \operatorname{Tr}(A):=\{F\in A^*\,|\,F(ab)=F(ba)~~\forall~a,b\in A\}.
\]
Then for $F,G\in\operatorname{Tr}(A)$ it is a well-known consequence of
the fact that $\Delta$ is a morphism of associative algebras
that the result of the operation
\begin{equation}\label{eqn:QuasiDeltaAssOnTraces}
   F\bullet G:=(F\otimes G) \Delta
\end{equation}
is again a trace, and eqn (\ref{eqn:coass}) implies
that $(\operatorname{Tr}(A),\bullet,\epsilon)$ is an \textit{associative
algebra}.

\subsection{Quasi-Hopf quantized universal envelopping algebras}
 QLBA are classical objects of quasi-Hopf quantized universal envelopping
 algebras (QHQUE algebra in the sequel) which are formal deformations of
 universal envelopping algebras in the category of QH algebras. More precisely
 a QHQUE algebra is a QH algebra $A_\h = (A,\mu,\Delta,\Phi)$ over the ring of formal
 power series $\sfo\K$ (where tensor products are completed in the $\h$-adic topology)
 such that $A \simeq \sfo{(A/\h A)}$,
 $\Phi \equiv 1\otimes 1\otimes 1 \mod\h^2$ and:
\begin{equation*}
 A_\h/\h A_\h \simeq (U(\ggo),\mu_0,\Delta_0,1\otimes 1\otimes 1)
\end{equation*}
as QH algebras. Under these conditions (where $(a\otimes b)^{21}:=b\otimes
a$)
\begin{equation}
 \delta(x) \equiv (\Delta(x)-\Delta(x)^{21})/\h \mod\h\et\varphi
 \equiv (\alt\Phi)/\h^2 \mod\h \label{eqn:cl_1}
\end{equation}
with $x\in\ggo$, define a structure of QLBA on $\ggo$ (see \cite{Dri90} Proposition
2.1).\\
In the important particular case where the algebra $U(\ggo)$ is \textit{rigid}, i.e. the
multiplication of $U(\ggo)$ can be regarded as undeformed in $A_\h$, it
thus follows that the space of traces of $A_\h$, $\operatorname{Tr}(A_h)$, is equal to
$\operatorname{Tr}\big(U(\ggo)\big)[[\h]]$ as a $\sfo\K$-module. The results of
the preceeding subsection imply that in this particular case the
trace-algebra
$(\operatorname{Tr}(A_h),\bullet,\epsilon)$ of the QHQUE associated to the
QLBA $(\ggo,[\,,],\delta,\varphi)$ is a
\textit{deformation quantization of the associated Poisson algebra}
$\big(\operatorname{Tr}(U(\ggo)),\bullet_0,\{\,,\}_D\big)$.

\subsection{Twist transformation}
When talking about equivalence between QLBA or between QH algebras,
the proper notion is the twist transformation. This can be used to
obtain various quantizations of QLBA by twisting QHQUE.
Let $\ggo=(\ggo,\delta,\varphi)$ be a QLBA, the twist of $\ggo$
via a skew element $f\in\bw^2 \ggo$ is the QLBA $\ggo^f=(\ggo,\delta^f,\varphi^f)$ where:
\begin{equation}
\delta^f(x) = \delta(x)+[x\otimes 1+1\otimes x,f]\et
\varphi^f = \varphi+\cp(\delta\otimes\id)(f)-\cyb(f)\label{eqn:twist_lqba}
\end{equation}
with $x\in\ggo$ and $\cyb(f)=[f^{12},f^{13}]+[f^{12},f^{23}]+[f^{13},f^{23}]$.
Twisting via $f_1+f_2$ is equivalent to twisting first via $f_1$, then via $f_2$.
Let $(A,\mu,\Delta,\Phi)$ be a QH algebra, the twist of $A$ via an invertible
element $F\in A^{\otimes 2}$ is the QH algebra $A^F=(A,\mu,\Delta^F,\Phi^F)$ where:
\begin{equation*}
 \Delta^F(x) = F\Delta(x)F^{-1}\et
\Phi^F = F^{23}F^{1,23}\Phi (F^{-1})^{12,3}(F^{-1})^{12}
\end{equation*}
with $x\in A$ and $F^{12,3}$, $F^{1,23}$ have the obvious meaning $(\Delta\otimes\id)(F)$,
$(\id\otimes\Delta)(F)$. Twisting via $F_2 F_1$ is equivalent to twisting first via
$F_1$, then via $F_2$.\\
If $A_\h$ is a QHQUE algebra with classical limit $\ggo$, then the QHQUE algebra
$A_\h^{F}$ obtained by twisting $A_\h$ via $F\in\sfo{U(\ggo)^{\otimes 2}}$ is
$\ggo^f$ where:
\begin{equation*}
 f \equiv (F^{21}-F)/h \mod h.
\end{equation*}
\begin{rmk}
 The only non-trivial point is to show that the element $f$ below, which
 \emph{a priori} lies in $U(\ggo)^{\otimes 2}$, actually belongs to $\bw^2\ggo$.
 This is proved by using a  cohomological argument, see \cite{Dri90}, Proposition 2.1.
\end{rmk}

\section{Pohlmeyer-Rehren QLBA and their quantization}\label{sec:3}

Let $V$ be a finite dimensional $\K$-vector space, $T(V)$ the tensor algebra over
$V$ and $L(V)$ the free Lie algebra over $V$. We recall that $U(L(V))=T(V)$ and
that $L(V)$ is exactly the set of primitive elements of $T(V)$, \emph{i.e.}
\begin{equation*}
 x\in T(V)\ \textnormal{s.t.}\ \Delta_0(x) = x\otimes 1+1\otimes x.
\end{equation*}

\subsection{Prerequisites}\label{sec:der}
Let us recall some facts about derivations over a morphism. Let $A$, $A'$, $B$, $B'$ and $C$ be
algebras. Let $D\in\Hom_\K(A,B)$ and $\phi\in\Hom_{\operatorname{Alg}}(A,B)$. We say
that  $D$ is a derivation over $\phi$ if $\forall x,y\in A$:
\begin{equation}
 D(xy) = D(x)\phi(y)+\phi(x)D(y).\label{eqn:der}
\end{equation}
We denote by $\Hom_\phi(A,B)$ the set of derivations over $\phi$. There are several
(and straightforward to prove) properties of such derivations:
\begin{align}
 \left.\begin{array}{c}
            D\in\Hom_\phi(A,B)\\
            \psi\in\Hom_{\operatorname{Alg}}(B,C)
           \end{array}\right\rbrace & \Longrightarrow \psi D\in\Hom_{\psi\phi}(A,C),
             \label{eqn:der_2}\\
 \left.\begin{array}{c}
      D\in\Hom_\phi(A,B)\\
      \psi\in\Hom_{\operatorname{Alg}}(A',B')
     \end{array}\right\rbrace & \Longrightarrow D\otimes\psi\in\Hom_{\phi\otimes\psi}
       (A\otimes A', B\otimes B'),\label{eqn:der_3}\\
 D_1,D_2\in\Hom_\phi(A,B) & \Longrightarrow D_1+D_2\in\Hom_\phi(A,B).\label{eqn:der_4}
\end{align}
Properties \refeqn{eqn:der_2} and \refeqn{eqn:der_3} also work by composing/tensoring
in the other way around. It is well known that for any $\K$-morphisms $f$ from $V$ to an
associative algebra $A$, there exists an unique morphism of algebras $F$ from $T(V)$
to $A$ which coincide with $f$ on $V$. We have a similar property for derivations over
a morphism.

\begin{lem}\label{lemma:up_der}
 For $d,\phi\in\Hom_\K(V,A)$ their exists an unique $D\in\Hom_{\Phi}(T(V),A)$ such
 $D|_V=d$ and $\Phi|_V=\phi$.
\end{lem}
\bproof
We first extend $\phi$ to $\Phi$ using the universal property of the tensor algebra
$T(V)$, then we use \refeqn{eqn:der} as extension formula, and since it respects the
associativity of the product of $T(V)$, the lemma is proved.
\eproof
\begin{rmk}
 This lemma will be useful in the sequel: in order to prove formulas on $T(V)$
 where both sides are derivations over the same morphism, it is sufficient to check the
 formulas on the generators $v\in V$.
\end{rmk}

\subsection{Pohlmeyer-Rehren quasi-Lie bialgebras}

\begin{thm}
 Let $s\in V^{\otimes 2}$, then the formulas:
\begin{equation}
 \delta_s(x) = [s,x\otimes 1]-[s^{21},1\otimes x]\et\varphi_s = -\cp[s^{12},s^{13}]
  \label{eqn:def_lqba_1}
\end{equation}
with $x\in V$, define a QLBA structure on $L(V)$. Moreover, $\delta_s$ satisfies
co-Jacobi if, and only if $s=v\otimes w$ for some $v,w\in V$.
\end{thm}
\bproof
It is clear that $\delta(x)\in \bw^2 L(V)$. We use Lemma \ref{lemma:up_der} to
extend $\delta$ to a derivation $D:T(V)\longrightarrow T(V)^{\otimes 2}$ over the
shuffle comultiplication $\Delta_0$  on $T(V)$ (generated by the diagonal map
$x\longmapsto x\otimes 1+1\otimes x$ on $V$). It is not clear that $D$ can be
restricted to a map $D:L(V)\rightarrow L(V)^{\otimes 2}$. To prove this we have
to show that $D$ satisfies the co-Leibniz rule on $T(V)$ (see \refeqn{eqn:coleibniz}
for notations):
\begin{equation}
 (\Delta_0\otimes\id)D = (\id\otimes D)\Delta_0+\sigma_{23}(D\otimes\id)\Delta_0.
 \label{eqn:coleibniz_2}
\end{equation}
The left hand side of \refeqn{eqn:coleibniz_2} is a derivation over
$(\Delta_0\otimes\id)\Delta_0$ by \refeqn{eqn:der_3} and \refeqn{eqn:der_2}.
The right hand side is also a derivation over $(\Delta_0\otimes\id)\Delta_0$
by \refeqn{eqn:der_2}-\refeqn{eqn:der_4} (we recall that $\Delta_0$ is cocommutative
and coassociative). It is sufficient to check \refeqn{eqn:coleibniz_2} for $x\in V$:
\begin{align}
 (\Delta_0\otimes\id)D(x) &= [s^{13},x_1]+[s^{23},x_2]-[s^{31}+s^{32},x_3]
       \label{eqn:col_proof_1}\\
(\id\otimes D)\Delta_0(x) &= 1\otimes D(x)
=[s^{23},x_2]-[s^{32},x_3]\label{eqn:col_proof_2}\\
\sigma_{23}(D\otimes\id)\Delta_0(x) &= \sigma_{23}(D(x)\otimes 1)
= [s^{13},x_1]-[s^{31},x_3],\label{eqn:col_proof_3}
\end{align}
and we have $\refeqn{eqn:col_proof_1}=\refeqn{eqn:col_proof_2}+\refeqn{eqn:col_proof_3}$.
Let $x\in L(V)$, we need to prove that $D(x)$ lies in $L(V)^{\otimes 2}$. We can write
$D(x)$ as a sum of tensors $\sum_i x'_i\otimes x''_i$ such that the $x''_i$'s are free,
applying the co-Leibniz rule gives:
\begin{equation*}
 \sum_i \Delta_0(x'_i)\otimes x''_i = (x'_i\otimes 1+1\otimes x'_i)\otimes x''_i.
\end{equation*}
It follows that the $x'_i$'s are primitives, \emph{i.e.} lies in $L(V)$. Since $D$ is skew:
\begin{equation*}
 D(L(V))\subset (L(V)\otimes T(V))\cap (T(V)\otimes L(V))=L(V)\otimes L(V).
\end{equation*}
We call $\delta$ (again) the restriction of $D$ to $L(V)$, and it is clear from
\refeqn{eqn:der} that $\delta\in Z^1(L(V),\bw^2 L(V))$. Now we want to show that
$\delta$ satisfies \refeqn{eqn:cojac}, for this purpose we will show:
\begin{equation}
 \cp(D\otimes\id)D(x) = -\big[(\Delta_0\otimes\id)\Delta_0(x),\cp[s^{12},s^{13}]\big]
 \label{eqn:cojac_proof_1}
\end{equation}
with $x\in T(V)$. Again the two sides of \refeqn{eqn:cojac_proof_1} are derivations
over $(\Delta_0\otimes\id)\Delta_0$ (using \refeqn{eqn:der_2}-\refeqn{eqn:der_4}).
On the generators $x\in V$, \refeqn{eqn:cojac_proof_1} follows from a straightforward
computation. The last condition $\alt(\delta\otimes\id\otimes\id)(\cp[s^{12},s^{13}])=0$
is a long and straightforward computation on permutations.
\eproof

\begin{defi}\label{def:pr_lqba}
 The Pohlmeyer-Rehren QLBA is $(L(V),\delta_g,\varphi_g)$ with $g\in S^2(V)$ of
 signature $(-1,1,..,1)$, namely:
\begin{equation*}
 \delta_g(x)=[g,x\otimes 1-1\otimes x]\et
 \varphi_g = -[g^{12},g^{13}]+[g^{12},g^{23}]-[g^{13},g^{23}].
\end{equation*}
\end{defi}

To quantize $L(V)$ we need to deform $U(L(V))=(T(V),\mu_0,\Delta_0,1\otimes 1\otimes 1)$
in the category of QH algebras. It is well known that we do not need to deform the unit
$1$ nor the counit $\eps_0$ (see \cite{SS93} for example). Moreover the multiplication
$\mu_0$ is rigid. The only data we have to take care of are the coproduct and the element
$\Phi$.

\paragraph{Coboundary case.} As a first step towards the quantization, it is interesting to
see what happens in the coboundary case. Choose $s\in \bw^2 V$, since $s^{21}=-s$:
\begin{equation*}
 \delta_s(x)=[s,x\otimes 1+1\otimes x]\et
\varphi_s = -\cp[s^{12},s^{13}] = -\cyb(s).
\end{equation*}
Consider $(L(V),\delta_s^s,\varphi_s^s)$ the twist of $(L(V),\delta_s,\varphi_s)$ via $s$
(see \refeqn{eqn:twist_lqba}), we have:
\begin{align*}
& \delta_s^{s}(x) = [s,x\otimes 1+1\otimes x]-[s,x\otimes 1+1\otimes x] = 0,\\
& \varphi_s^{s} = -\cyb(s)+\cp(\delta\otimes\id)(s)-\cyb(s) = -2 \cyb(s)+2\cyb(s) = 0.
\end{align*}
Hence $(L(V),\delta_s,\varphi_s)$ can be twisted to the trivial QLBA structure on $L(V)$.
 A quantization of $(L(V),0,0)$ would be
 $A'_\h:=(\sfo{T(V)},\mu_0,\Delta_0,1\otimes 1\otimes 1)$ and consequently a
 quantization $A_\h$ of $(L(V),\delta_s,\varphi_s)$ would be the twist of $A'_\h$ via
\begin{equation*}
 F=e^{\h \frac{s}{2}}\in \sfo{T(V)^{\otimes 2}}.
\end{equation*}
 Namely $A_\h:=(\sfo{T(V)},\mu_0,\Delta,\Phi)$:
\begin{equation*}
 \Delta(x) = F(x\otimes 1+1\otimes x)F^{-1}\et \Phi =
 F^{23}F^{1,23}(F^{-1})^{12,3}(F^{-1})^{13}
\end{equation*}
with $x\in T(V)$. Quantization of coboundary QLBA is already known
(see \cite{EH05} and \cite{EH06}), however we have here a simple and
explicit quantization.

\subsection{Low rank cases}

In this subsection we study what happens if the bivector $g$ is of low rank $\leq 2$,
which is unphysical since $g$ comes from a
Minkowski metric (see section \refeqn{sec:1}) in an at least $3$-dimensional
vector space.
We set $s=v\otimes w$, $g=(s+s^{21})/2$ and $f=(s-s^{21})/2$ with $v,w\in V$. It is
clear that:
\begin{equation}
  [s^{12},s^{13}] = [s^{13},s^{23}] = 0\et [s^{12},s^{23}] \neq 0.\label{eqn:rk2_1}
\end{equation}
From \refeqn{eqn:def_lqba_1} we see that $(L(V),\delta_s)$ is a genuine Lie bialgebra.

\begin{prop}
 An explicit quantization of $(L(V),\delta_s)$ is the Hopf algebra
 $A'_\h=(\sfo{T(V)},\mu_0,\Delta',1\otimes 1\otimes 1)$:
\begin{equation*}
 \Delta'(x) = G(x\otimes 1)G^{-1}+1\otimes x,
\end{equation*}
with $x\in V$ and $G=e^{\h s}\in\sfo{T(V)^{\otimes 2}}$.
\end{prop}

\bproof
 From \refeqn{eqn:rk2_1} we know that the only non-commuting terms in $G^{12},G^{13}$
 and $G^{23}$ are $G^{12}$ and $G^{23}$. Moreover $v$ is primitive for $\Delta'$ but
 $w$ is not (since $[s,w\otimes 1]\neq 0$). Let start with computing $G^{1,23}$ and
 $G^{12,3}$:
\begin{align}
& G^{1,23} = e^{\h (v\otimes\Delta'(w))} = e^{\h (G^{23} s^{12} (G^{-1})^{23} + s^{13})}
 = G^{23}G^{12}(G^{-1})^{23}G^{13},\\
& G^{12,3} = e^{\h (\Delta'(v)\otimes w)} = e^{\h (s^{13}+s^{23})} = G^{13}G^{23}.
\end{align}
Then:
\begin{align}
 (\id\otimes\Delta')\Delta'(x) =& G^{13}G^{23}G^{12} x_1
 (G^{-1})^{12}(G^{-1})^{23}(G^{-1})^{13}\label{eqn:rk2_3}\\
 & +G^{13}G^{23}x_2 (G^{-1})^{23}(G^{-1})^{13}+x_3,\nonumber\\
 (\Delta'\otimes\id)\Delta'(x) =& G^{23}G^{12}(G^{-1})^{23}G^{13} x_1
 (G^{-1})^{13}G^{23}(G^{-1})^{12}(G^{-1})^{23}\label{eqn:rk2_4}\\
 & +G^{23}x_2 (G^{-1})^{23}+x_3.\nonumber
\end{align}
To complete the proof we just have to notice that the non-commuting terms $G^{12}$
and $G^{23}$ (and their inverse) are in the same order in \refeqn{eqn:rk2_3} and
in \refeqn{eqn:rk2_4}.
\eproof

$(L(V),\delta_g,\varphi_g)$ is the twist of the QLBA $(L(V),\delta_s,\varphi_s)$ via
the skew element $f$. Indeed:
\begin{align*}
 & \delta_s^f(x) = [s,x_1]-[s^{21},x_2]+[x_1+x_2,(s-s^{21})/2] = [(s+s^{21})/2,x_1-x_2],\\
 & \varphi_s^f = \cp(\delta_s\otimes\id)(f)-\cyb(f) = 2\cyb(f)-\cyb(f) = \varphi_g.
\end{align*}
It follows that a quantization $A_\h$ of $(L(V),\delta_g,\varphi_g)$ would be the twist
of $A'_\h$ via:
\begin{equation*}
 J = e^{-\h \frac{s}{2}} = \sqrt{G^{-1}}\in\sfo{T(V)^{\otimes 2}}.
\end{equation*}
Namely $A_\h:=(\sfo{T(V)},\mu_0,\Delta,\Phi)$:
\begin{equation*}
  \Delta(x) = J^{-1}(x\otimes 1)J+J(1\otimes x)J^{-1}\et
  \Phi = (J^{-1})^{23}J^{12}J^{23}(J^{-1})^{12}
\end{equation*}
with $x\in V$. We just proved the
\begin{thm}
 $A_\h$ is a quantization of $(L(V),\delta_g,\varphi_g)$.
\end{thm}

\begin{rmk}
 The particular case $w=\lambda v$ for some $\lambda\in k$ yields the case
 where $g$ has rank 1.
\end{rmk}

Moreover, we compute the deformed antipode for $A'_\h$.
The quasi-antipode $(S,\alpha,\beta)$ for $A_\h$ can be obtained by using the
twisting formulas for QH algebras (see \cite{Dri90}).
\begin{prop}
 The antihomomorphism $S:\sfo{T(V)}\rightarrow\sfo{T(V)}$ defined on generators $x\in V$
 by:
\begin{equation}
 S(x) = -e^{-\h \Lm_v\Rm_w}(x)\gamma^{-1}\label{antipode}
\end{equation}
(where $\Lm_a(b):=ab=:\Rm_b(a)$ for all $a,b\in T(V)$ and
$\gamma := \mu_0(G^{-1})$) is an antipode for $A'_\h$.
\end{prop}
\bproof
 The antipode of $T(V)$ is $S_0(x)=-x$ (and extended as antihomomorphism) for $x\in V$.
 Consequently the deformed antipode $S$ exists and is unique. The convolution product
 $\star$ on $\operatorname{End}_{\sfo\K}(A_\h)\simeq \sfo{\operatorname{End}_\K(T(V))}$
 is associative and $1\eps_0$ is an unit for $\star$. Then we just have to find $S$
 satisfying $S\star\id=1\eps_0$. First notice that $S(v)=-v$ since $v$ is a primitive
 element. Let $x\in V$:
\begin{equation*}
 (S\star\id)(x) = \sum_{r\geq 0}\h^r \sum_{i+j=r}\frac{(-1)^i}{i!j!} v^j S(x)
 v^i w^{i+j}+x \overset{!}{=} 0.
\end{equation*}
Isolating $S(x)$ gives (\ref{antipode}). One still has to check if $S\star\id=1\eps$
which is a straight-forward computation.
\eproof

\subsection{General case}

The previous subsection shows how the difficulty is growing with the rank of the
bivector $g$. In the following, we can see how difficult can be the combinatorics of
rearranging legs of successive powers of $g$ order by order.

\begin{thm}
 Let $(L(V),\delta_g,\varphi_g)$ be a Pohlmeyer-Rehren QLBA (as in definition
 \ref{def:pr_lqba}).
 
 Then $(\sfo{T(V)},\mu_0,\Delta,\Phi)$ is a quantization
 of $L(V)$ modulo $\h^3$:
\begin{equation*}
 \Delta = \Delta_0+\h \Delta_1+\h^2 \Delta_2+O(\h^3)\et
  \Phi = 1\otimes 1\otimes 1-\frac{1}{2} h^2 [g^{12},g^{13}]+O(h^3)
\end{equation*}
where:
\begin{equation*}
 \Delta_0(x) = x\otimes 1+1\otimes x,\quad \Delta_1(x) = [g,x\otimes 1]
\end{equation*}
\begin{align*}
 \Delta_2(x) = & \alpha g^2 x_1 - \frac{1}{2}\ g x_1 g + \left(\frac{1}{2}-\alpha\right) x_1 g^2\\
 & + \left(\frac{1}{2}-\alpha\right)(\id\otimes\tau)(g^2 x_1) - \frac{1}{2} (\id\otimes\tau)(g x_1 g) +
 \alpha (\id\otimes\tau)(x_1 g^2)\\
 & +\left(\beta-\frac{1}{2}\right) g^2 x_2 + \frac{1}{2}\ g x_2 g - \beta\ x_2 g^2\\
 & + \left(\frac{1}{2}-\beta\right)(\tau\otimes\id)(g^2 x_2)-\frac{1}{2}(\tau\otimes\id)(g x_2 g)+
 \beta (\tau\otimes\id)(x_2 g^2)
\end{align*}
with $x\in V$, $\alpha,\beta\in\K$ and $\tau(x_1 x_2)=x_2 x_1$.
\end{thm}

\bproof
We only sketch the proof which is a very long and tedious computation.
The first step is the ansatz that $\Delta_2(x)$ is a linear combination of
the 12 terms $g^2 x_1$, $g^2 x_2$, $(\id\otimes\tau)(g^2 x_1)$,
$(\tau\otimes\id)(g^2 x_2)$,\ldots. Then we compute \refeqn{eqn:coass}
at order 2 and find a linear system with constant terms depending on
12 parameters and several equations. The rank of this system is 10, and
the number of equations can be reduced to 10. Consequently the
space of solutions is of dimension $2$.
\eproof

\section{Invariant charges by Pohlmeyer-Rehren}\label{sec:1}

The goal of this section is to roughly present the theory of reparametrization-invariant charges due to Pohlmeyer-Rehren and the Poisson brackets on the $\zed_{\mu_1\cdots\mu_k}^\pm$ symbols. For further details see \cite{Poh86}, \cite{PR86}, \cite{Nam70} and \cite{Got71}.

According to the Nambu-Goto theory, the world sheet of a closed string moving in
a Minkowski space $(\mink,g=\langle\,,\rangle)$, is a smooth map
$x:\Rb\times\SXb^1\to\mink$ such that the induced metric $\til{g}[x]:=x^* g$
is of signature $(-1,1)$ ($d\geq 2$) and such that the area is extremal with
respect to $\til{g}[x]$:
\begin{equation}
 S[x] := M^2 \int_{\Rb\times\SXb^1} \sqrt{-\det\til{g}[x]_{(\tau,\sigma)}}
 \dd\tau \dd\sigma\label{string_equation}
\end{equation}
where $M$ is a normalizing constant. We set $\dot x := \dep{x}{\tau}$,
$x':=\dep{x}{\sigma}$ and we define:
\begin{equation*}
 u_\mu^\pm := p_\mu\pm M^2 \left\langle x',e_\mu\right\rangle\et
 p_\mu := \frac{M^2}{\sqrt{-\det\til{g}[x]}}\left(\left\langle x',
 \dot x\right\rangle \left\langle x',e_\mu\right\rangle-
 \left\langle x',x' \right\rangle \left\langle \dot x,e_\mu\right\rangle\right)
\end{equation*}
with
$\{e_\mu\}_{\mu=0}^{d-1}$ a basis of $\mink$. By making a zero-curvature ansatz
for the Euler-Lagrange equations derived from \refeqn{string_equation}
Pohlmeyer gets the $\er^\pm$ and $\zed^\pm$-quantities which are iterated
integrals of the following form:
\begin{align*}
& \er^\pm_{\mu_1\cdots\mu_k}(\tau,\sigma)
:= \int_{\sigma\leq\sigma_k\leq\cdots\leq\sigma_1\leq\sigma+2\pi}
 u^{\pm}_{\mu_1}(\tau,\sigma_1)\cdots
 u^{\pm}_{\mu_k}(\tau,\sigma_k)\dd\sigma_1\cdots\dd\sigma_k,\\
& \zed^\pm_{\mu_1\cdots\mu_k} := \er^\pm_{\mu_1\cdots\mu_k}
+\er^\pm_{\mu_2\cdots\mu_k\mu_1}+\cdots+\er^\pm_{\mu_k\mu_1\cdots\mu_{k-1}}.
\end{align*}
The cyclically invariant $\zed^\pm$-quantities do no longer depend on $\sigma$
and $\tau$ (and are thus
in particular conserved quantities for the moving string) and are shown to be
reparametrization invariant observables.

Pohlmeyer and Rehren \cite{PR86} have found the following Poisson structure
on the vector space
of all the $\zed^\pm$-observables: first the pointwise multiplication
of two iterated integrals is well-known to be an iterated integral
{\small\begin{equation*}
 \er^\pm_{\mu_1\ldots\mu_k}(\tau,\sigma)
 \er^\pm_{\mu_{k+1}\ldots\mu_{k+l}}(\tau,\sigma)
 =  \sum_{\pi\in\ush(k,l)} \er^\pm_{\mu_{\pi(1)}\ldots\mu_{\pi(k+l)}}(\tau,\sigma)
    =: \er_{\boite{\mu_1\ldots\mu_k}{\mu_{k+1}\ldots\mu_l}}(\tau,\sigma).
\end{equation*}}
Here $\ush(k,l)$ denotes the set of all inverses of $(k,l)$-shuffle
permutations,
and where the box notation is due to Pohlmeyer-Rehren
(see \cite{PR86} for a proof and \cite{Che61,Che77a,Che77b},
\cite{Aom78} for details on iterated integrals). The same equation holds
for the $\zed^\pm$-quantities. \\
Let $\mathcal{P} = C^{\infty}(\SXb^1,\mink\times\mink)$
be the phase space (\emph{i.e.} the space of initial conditions
$\sigma\mapsto \big(x(0,\sigma),p(0,\sigma)\big)$ of the strings). The
$\er^\pm$-quantities can be regarded as maps from $\mathcal{P}$
to $\mathcal{C}^\infty(\SXb^1,\mathbb{R})$, and the
$\zed^\pm$-quantities as maps from $\mathcal{P}$ to $\mathbb{R}$.
By regularizing the physical rule of thumb called \textit{fundamental equal time Poisson
brackets}
$\{x_\mu(\tau,\sigma),p_\nu(\tau,\sigma')\}$``=''$g_{\mu\nu}\delta(\sigma-\sigma')$
(mimicking a symplectic structure on $\mathcal{P}$)
Pohlmeyer and Rehren succeeded in finding the following Poisson bracket \cite{PR86}
for the $\zed^\pm$-observables which satisfies the Jacobi identity:
{\small
\begin{equation}\label{eqn:PoissonBracketPohlmeyerRehren}
\quad\left\lbrace \zed_{\mu_1\ldots\mu_k}, \zed_{\nu_1\ldots\nu_{l}}
\right\rbrace_{\operatorname{PR}}
= 2\sum_{i=1}^k \sum_{j=1}^{l} g_{\mu_i \nu_j}
\left(\zed_{\mu_{i+1}\boite{\mu_{i+2}\ldots \mu_{i-1}}{\nu_{j+1}
        \ldots\nu_{j-2}}\nu_{j-1}}\right.
       \left.-\zed_{\nu_{j+1}\boite{\mu_{i+1}\ldots
       \mu_{i-2}}{\nu_{j+2}\ldots\nu_{j-1}}\mu_{i-1}} \right)
\end{equation}}
The r.h.s. is also definable for the $\er^{\pm}$'s, but \textit{fails to satisfy the Jacobi
identity}, a phenomenon which is well-known in the classical field theory of
certain integrable models.

In order to make contact to the theory of QLBAs and QHs, let $V$ be the
dual space of $\mink$, i.e. $V:=(\mink)^*$, and consider the free algebra
$T(V)$ and the free Lie algebra $L(V)$. In the preceeding paragraph we have seen
that the Minkowski metric $g\in V\otimes V$ gives rise to the Pohlmeyer-Rehren
$1$-cocyle $\delta_g$, and that therefore $\big(L(V),[\,,],\delta_g,\varphi_g\big)$
is a QLBA. Now
the free algebra $T(\mink)$ over the Minkowski space $\mink$
is a subspace of the dual space $T(V)^*$. Define the space of
\textit{restricted traces} $\operatorname{Tc}\left(T(V)\right)$ by
$\operatorname{Tr}\left(T(V)\right)\cap T(\mink)$ which is clearly
isomorphic to the subspace of \textit{cyclic} tensors in $T(\mink)$,
i.e. those elements whose homogeneous components are invariant under
all cyclic permutations.

By Proposition \ref{PropPoissonOnTraces} it is not hard to see that
$\big(\operatorname{Tc}\left(T(V)\right),\bullet_0,\{\,,\}_D\big)$
is a Poisson algebra where $\{\,,\}_D$ denotes the Poisson bracket
(\ref{eqn:MultShuffleUGBracketD}) with respect to the derivation extending
$\delta_g$.

Let us go back to the $\er_{\mu_1\cdots\mu_k}^\pm$ and
$\zed_{\mu_1\cdots\mu_k}^\pm$-symbols: they are obviously
paramet\-rized by $T(\mink)$ and can thus be seen as linear maps
\begin{eqnarray*}
 \er^\pm : T(\mink) & \longrightarrow & \operatorname{Fun}
 \left(\mathcal{P},C^\infty(\SXb^1,\Rb)\right)\\
 \zed^\pm : \operatorname{Tc}\left(T(V)\right)
    & \longrightarrow & \operatorname{Fun}
 \left(\mathcal{P},\Rb\right)
\end{eqnarray*}
by setting
\begin{equation*}
 \langle \er^{\pm}, e_{\mu_1}\cdots e_{\mu_k} \rangle (x,p): =
 \er^{\pm}_{\mu_1\cdots\mu_k}\et
 \langle \zed^{\pm}, e_{\mu_1}\cdots e_{\mu_k} \rangle (x,p): =
 \zed^{\pm}_{\mu_1\cdots\mu_k}.
\end{equation*}
Then we have the following Proposition which links the Pohlmeyer-Rehren
Poisson structure with the purely combinatorial Poisson structure defined
by the PR-$1$-cocycle:
\begin{prop}\label{PPoRehQuasi}
 Let $D$ be the derivation extension of the PR-$1$-cocycle $\delta_g$
 (see Definition \ref{def:pr_lqba}). Let $a,b\in T(\mink)$ and
 $a',b'\in \operatorname{Tc}(T(V))$. Then
\begin{eqnarray*}
 & \langle\er^\pm,a\rangle\langle\er^{\pm},b\rangle = \langle \er^{\pm},a\bullet_0 b\rangle,\quad
 \langle\zed^\pm,a'\rangle\langle\zed^{\pm},b'\rangle = \langle \zed^{\pm},a'\bullet_0 b'\rangle, &
 \\
 &\et
  \left\lbrace \langle\zed^\pm,a'\rangle,\langle\zed^\pm,b'\rangle
  \right\rbrace_{\operatorname{PR}} = \left\langle
  \zed^\pm,\{a',b'\}_D\right\rangle. &
\end{eqnarray*}
 It follows that
 $\zed^\pm$ are Poisson maps from $\operatorname{Tc}(T(V))$ to its image.
\end{prop}
\bproof
Sketch: as it has been mentioned above, the link of the pointwise
multiplication with the unshuffle multiplication $\bullet_0$ is a
classical property of iterated integrals. The identity of the PR-Poisson
bracket (\ref{eqn:PoissonBracketPohlmeyerRehren}) with the algebraic
bracket $\{\,,\}_D$ is a straight-forward computation using the shuffle
comultiplication. The appearance of $g$ in the formula of $\delta_g$ in
Definition \ref{def:pr_lqba} is already visible in eqn
(\ref{eqn:PoissonBracketPohlmeyerRehren}) where $g_{\mu_i\nu_j}$ occurs
in front of the $\zed^{\pm}_{\mu_1\cdots\mu_k}$.
\eproof

As a possible framework for the deformation quantization of the Pohlmeyer-Rehren
Poisson algebra we thus propose the quantization of the Pohlmeyer-Rehren QLBA
as a QHQUE: according to eqn (\ref{eqn:QuasiDeltaAssOnTraces}) the restricted traces
will become an associative algebra parametrizing the
$\zed^\pm$-observables.

\end{document}